\documentclass[12pt]{article}
\usepackage[T2A]{fontenc}
\usepackage[russian]{babel}
\usepackage{amsfonts}
\usepackage{amsmath}
\usepackage{amssymb}
\usepackage{amsthm}

\renewcommand{\div}{\opn{div}}
\newcommand\pa{\partial}

\newcommand\ov{\overline}

\newcommand\opn{\operatorname}

\newcommand\D{{\mathcal D}}
\renewcommand\H{{\mathcal H}}

\newcommand\R{{\mathbb R}}
\newcommand\N{{\mathbb N}}
\newcommand\Z{{\mathbb Z}}

\newcommand\al{\alpha}
\newcommand\be{\beta}

\newcommand\De{\Delta}

\newcommand\la{\lambda}

\newcommand\om{\omega}
\newcommand\Om{\Omega}

\renewcommand\th{\theta}

\renewcommand\phi{\varphi}

\newcommand\ep{\varepsilon}

\renewcommand\({\left(}
\renewcommand\){\right)}

\def\H{{\mathfrak H}}

\def\suml{\sum\limits}

\newcommand{\const}{\opn{const}}

\newlength{\MMM}
\newlength{\MMMM}
\newcommand{\Ho}{%
\makebox{\parbox[l][\MMMM][b]{\MMM}{%
\raisebox{0.12em}{$\stackrel{\circ}{H}$}}}}

\renewcommand{\L}{\mathcal L}
\renewcommand{\N}{\mathcal N}

\newtheorem{theor}{Теорема}[section]
\newtheorem{theorem}[theor]{Теорема}
\newtheorem{proposition}[theor]{Предложение}
\newtheorem{corollary}[theor]{Следствие}
\newtheorem{lemm}[theor]{Лемма}

\newenvironment{proof*nodot}{\par\smallskip{\bf Доказательство}\rm}
			{\hfill $\square$\medskip}

\newenvironment{theorem*}[1]{\par\smallskip{\bf #1.}\it}{\par\medskip}

\begin{document}

%\large
\makeatletter
%\renewcommand{\@oddhead}{\hfil --- \thepage ---\hfil}
%\renewcommand{\@oddfoot}{}
%\@addtoreset{equation}{section}

%\renewcommand{\section}{\@startsection{section}{1}%
%{-1.4cm}{3.5ex plus 1ex minus .2ex}%
%{2.3ex plus .2ex}{\LARGE\bf\hskip0.8cm}
%}

\renewcommand{\subsection}{\@startsection{subsection}{2}%
{0pt}{2.25ex plus 1ex minus .2ex}%
{0.5ex plus .2ex}{\large\bf}
}

\makeatother

\renewcommand{\bibname}{Список литературы}

\title{Операторы Шредингера с сингулярными потенциалами в ограниченной
области}

\author{
М.~И.~Нейман--заде\footnote{Работа поддержана грантом
РФФИ N 01-01-00691}, А.~А.~Шкаликов\footnote{Работа поддержана грантом
РФФИ N 01-01-00958}}

\maketitle

Целью настоящей работы является обобщение результатов статьи~\cite{NeSh}
на случай ограниченной области $\Om\subset\R^n$ с гладкой границей.

Пусть $D(\Om)$ --- пространство тест--функций в $\Om$ (бесконечно гладких с
компактным носителем), а
$D'(\Om)$ --- соответствующее пространство распределений.
Мы ставим следующую задачу: для каких функций $q\in D'(\Om)$ корректно
определен оператор $-\De+q$?

Мы также будем изучать сильно эллиптические
операторы
\begin{equation}\label{III:eq:L}
L=\suml_{|\al|,|\be|\le m} D^\al c_{\al,\be}(x) D^\be,\qquad x\in \Om\subset \R^n,
\end{equation}
где коэффициенты $c_{\al,\be}(x)$ --- непрерывные функции при $|\al|=|\be|=m$, а
при $|\al|+|\be|<2m$, --- распределения.

Возникает аналогичный вопрос:
для каких функций $c_{\al,\be}\in D'(\Om)$
корректно определен оператор вида~\eqref{III:eq:L}? Здесь решение поставленных
задач осложняется тем, что в случае $\Om\ne\R^n$ для задания оператора нужно
задавать краевые условия. Самые простые краевые условия в нашей ситуации ---
условия Дирихле. Для условий Дирихле можно без труда получить аналоги
результатов обеих работ~\cite{NeSh},~\cite{Ne}. С другими краевыми условиями дело
обстоит сложнее. Здесь мы подробно разберем только случай оператора $-\De+q$
с обобщенными условиями Неймана.

В конце работы для операторов $-\De+q$, $q\in D'$, с условиями Дирихле и Неймана,
мы изучаем асимптотику собственных значений и базисные свойства собственных функций.

\subsection{Задача Дирихле для оператора Лапласа и сильно эллиптического
оператора, возмущенных сингулярными коэффициентами}

Сначала рассмотрим оператор $-\De+q$ и свяжем с ним квадратичную форму
$$
 ((-\De+q)\phi,\phi)=(\nabla\phi,\nabla\phi)+(q\phi,\phi),
$$
определенную на тест--функциях $\phi\in D(\Om)$. Очевидно, форма
$(\nabla\phi,\nabla\phi)$ допускает замыкание на пространство $\Ho^1(\Om)$.
Поэтому с выражением $-\De+q$ можно связать оператор, если
выполнено неравенство
$$
 |(q\phi,\phi)|\le\ep\|\phi\|_{H^1}+M\|\phi\|_{L_2},\quad\phi\in D,
$$
где $\ep$ достаточно мало, а $M=M(\ep)$. Мы должны замкнуть это
неравенство на функции $\phi\in\Ho^1(\Om)$. Но тогда это неравенство
эквивалентно тому, что $q$ есть мультипликатор из $\Ho^1(\Om)$ в дуальное
пространство по отношению к $L_2(\Om)$, причем относительная норма этого
мультипликатора достаточно мала. Дуальное  к $\Ho^1(\Om)$ пространство
обозначается через $H^{-1}(\Om)$ и для него справедлива
``структурная теорема'' (см.~\cite[гл.~1.12]{LM}):

\noindent{\it
Всякий элемент $f\in H^{-1}(\Om)$ представим в виде
\begin{equation}\label{IV:eq:struct}
f=\suml_{i=1}^n D_{x_i}f_i,
\end{equation}
где $f_i\in L_2(\Om)$, а равенство понимается в смысле
распределений. Норма в $H^{-1}$ эквивалентна норме
$$
 \|f\|^2_{-1}=\inf\suml_{i=1}^n\|f_i\|^2,
$$
где $\inf$ берется по всем $f_i$, для которых справедливо
представление~\eqref{IV:eq:struct}.
}

Далее, пространства $\Ho^s(\Om)$ при $0\le s\le 1$ можно определить
с помощью интерполяции $\Ho^s(\Om)=[\Ho^1(\Om),L_2(\Om)]_s$, и затем
положить $H^{-s}(\Om)=(\Ho^s(\Om))'$. Впрочем, существуют и другие способы
определения пространств $\Ho^s(\Om)$ и $H^{-s}(\Om)$ при $0\le s\le1$,
а также при всех $s\in\R$. С ними можно ознакомиться в
статьях Хермандера~\cite{Her1}, Волевича и Панеяха~\cite{VP}. Другие
определения эквивалентны указанному при $s\ne\pm1/2$
(в общем случае при $s\ne k+1/2$, $k\in\Z$).

Согласно теореме о следах пространства $\Ho^m_p(\Om)$ при целых $m\ge1$
корректно определены при $p>1$. Эти пространства состоят из функций
$u\in H_p^m(\Om)$, для которых
$$
 u(x)|_{\pa\Om}=\dfrac{\pa}{\pa\vec n}u(x)|_{\pa\Om}=\dots=
 \dfrac{\pa^{m-1}}{\pa\vec n^{m-1}}u(x)|_{\pa\Om}=0,
$$
где $\vec n$ есть нормаль к $\pa\Om$.
В книге~\cite{MSh} пространство $H^{-m}_p(\Om)$ определено с помощью
``структурной теоремы''. А именно, $H^{-m}_p(\Om)$ состоит из распределений
$u(x)\in D'(\Om)$, для которых справедливо представление
\begin{equation}\label{IV:eq:structm}
u(x)=\suml_{|\al|\le m} D^\al V_{\al}(x),\qquad
\text{где $V_\al\in L_p(\R^n)$}.
\end{equation}
Норму в $H^{-m}_p(\Om)$ можно определить равенством
$$
\|u\|_{-m,p}=\inf\(\suml_{|\al|\le m} \|V_\al(x)\|^p \)^{1/p},
$$
где $\inf$ берется по всем $V_\al$, для которых справедливо
представление~\eqref{IV:eq:structm}. При $\Om=\R^n$ такое определение
совпадает с прежним. При нецелых $s\le m$ пространства $H_p^{-s}$ можно
определить с помощью интерполяции.

Теперь вспомним доказательства теорем в статье~\cite{NeSh}.
Что было использовано в доказательствах, когда $\Om=\R^n$? Прежде всего,
определение
пространств $H_p^s(\R^n)$ и нормы в них. Теперь мы имеем такие определения
для $\Ho^s_p(\Om)$ и $H^{-s}_p(\Om)$. Далее, использовались
теорема вложения Соболева и некоторые результаты из книги
Мазьи--Шапошниковой (см. теоремы 6 и 7 статьи~\cite{NeSh}).
Но известно (см.~\cite{Maz},~\cite{Tr}), что теоремы
вложения сохраняются для $\Ho^s(\Om)$ и $H^{-s}(\Om)$, а
использованный результат Мазьи--Шапошниковой также справедлив
для гладкой области $\Om$. Следовательно, для мультипликаторов из
$\Ho^1(\Om)$ в $H^{-1}(\Om)$, а также из $\Ho^{m-|\al|}(\Om)$ в
$H^{-m+|\be|}(\Om)$, мы имеем такие же результаты, как при $\Om=\R^n$.
Сформулируем эти результаты.
\begin{theorem}\label{IV:th:mainDir2}
Пусть $\Om$ --- гладкая ограниченная область в $\R^n$. Пусть
\begin{equation}\label{IV:eq:qin}
q(x)\in
\left\{
\begin{aligned}
H^{-1}_2(\Om)&,\quad\text{если $n=1$},\\
H^{-1}_{2+\ep}(\Om)&,\quad\text{при некотором $\ep>0$, если $n=2$},\\
H^{-1}_n(\Om)&,\quad\text{если $n\ge 3$}.
\end{aligned}
\right.
\end{equation}
Тогда оператор $-\De+q(x)$ корректно определен с помощью метода квадратичных
форм. А именно, квадратичная форма
\begin{equation}\label{IV:eq:formL2}
(Lu,u)=(\nabla u,\nabla u)+(q(x), u\ov u),
\end{equation}
определенная при $u\in D(\Om)$, допускает замыкание на $\Ho^{1}(\Om)$ как
секториальная форма (в случае вещественной $q(x)$ --- как полуограниченная
форма). Следовательно, существует $m$--сек\-то\-риаль\-ный (полуограниченный)
оператор $L$ в $\H=L_2(\Om)$, ассоциированный с
этой формой. Этот оператор можно
также рассматривать как оператор в пространстве $\H^{-1}=H^{-1}(\Om)$
с областью определения $\D(L)=\H^1=\Ho^1(\Om)$. С этой точки зрения
$L$ есть компактное возмущение оператора $L_0=-\De$ в $\H^{-1}$
с той же областью $D(L_0)=\H^1$, причем $L$ имеет дискретный спектр.

Если $q_k(x)$ --- последовательность гладких функций таких, что при
$k\to\infty$
$$
\begin{aligned}
&\|q_k(x)-q(x)\|_{-1,2}\to 0,\quad\text{если $n=1$},\\
&\|q_k(x)-q(x)\|_{-1,2+\ep}\to 0,\quad\text{при некотором $\ep>0$, если
$n=2$},\\
&\|q_k(x)-q(x)\|_{-1,n}\to 0,\quad\text{если $n\ge 3$},
\end{aligned}
$$
то последовательность обычных операторов $L_k=-\De+q_k(x)$ в
$\H=L_2(\Om)$, порожденных условиями Дирихле, сходится в смысле
равномерной резольвентной сходимости к оператору $L$.
Собственные значения $\la_{s,k}$ операторов $L_k$ также сходятся к
собственным значениям $\la_s$ оператора $L$, причем скорость сходимости
оценивается через $\|q_k-q\|_{-1,p}$, где $p=2,2+\ep,n$ при $n=1,2$ и $\ge3$
соответственно.
\end{theorem}
\begin{proof}
Эта теорема есть следствие полученных нами ранее общих результатов об
операторах и теорем о мультипликаторах. Действительно, полагаем
$T=-\De$ в пространстве $\H^{-1}$ с областью определения $\D(T)=\H^1$.
Если выполнено условие~\eqref{IV:eq:qin} теоремы, то из теорем~6 и 7
статьи~\cite{NeSh} (точнее, аналогов этих теорем для $H^s(\Om)$) получаем,
что $q(x)$ компактный мультипликатор из $\Ho^1(\Om)=\H^1$ в
$H^{-1}(\Om)=\H^{-1}$. Согласно теореме~5 статьи~\cite{Ne}, оператор $T+q$
является $m$--секториальным в $\H^{-1}$ с той же областью, а его
квадратичная форма совпадает с замыканием формы~\eqref{IV:eq:formL2}.
Согласно результатам статьи~\cite{NeSh}, имеем вложение с оценкой норм
$$
 H^{-1}_p(\Om)\subset M[\Ho^1(\Om), H^{-1}(\Om)],
$$
где $p=2,2+\ep,n$ при $n=1,2$ и $\ge3$ соответственно. Поэтому равномерная
резольвентная сходимость операторов $L_n$ с гладкими потенциалами
к оператору $L$ и сходимость собственных значений следуют из
теоремы~3 статьи~\cite{NeSh} (если $q(x)$ вещественнозначна)
или из теоремы~5 статьи~\cite{Ne} (если $q(x)$ комплекснозначна).
Теорема доказана.
\end{proof}

Обобщение этой теоремы справедливо и для сильно эллиптических операторов
в $L_2$.

\begin{theorem}\label{IV:th:mainDirm}
Пусть оператор $L$ определен равенством~\eqref{III:eq:L}, где
$L_0$ --- равномерно сильно эллиптический оператор с непрерывными
коэффициентами $c_{\al,\be}(x)$ в области $\Om$. Пусть при $|\al|+|\be|<2m$
коэффициенты $c_{\al,\be}(x)$ таковы, что
\begin{align*}
 &c_{\al,\be}\in H^{|\al|-m}_p(\Om)\quad
\text{при $|\al|\ge|\be|$, $p>\max\{2,\dfrac{n}{m-|\be|}\}$},\\
 &c_{\al,\be}\in H^{|\be|-m}_p(\Om)\quad
\text{при $|\al|\le|\be|$, $p>\max\{2,\dfrac{n}{m-|\al|}\}$}.
\end{align*}
Тогда оператор $L$, определенный на тест--функциях $u(x)\in D(\Om)$,
допускает замыкание как секториальный оператор из пространства
$\Ho^1(\Om)=\H^1$ в пространство
$H^{-1}(\Om)=\H^{-1}$ (или как неограниченный оператор
в $\H^{-1}$ с областью определения $\D(L)=\H^1$). Оператор $L$ является
компактным возмущением оператора $L_0$, действующего в том же пространстве.
Если последовательности гладких функций $c_{\al,\be}^k(x)$ сходятся к
$c_{\al,\be}$ в соответствующих пространствах $H^{|\al|-m}_p(\Om)$,
то регулярные операторы $L_k$, порожденные краевыми условиями Дирихле,
сходятся к оператору $L$.
\end{theorem}
\begin{proof*nodot} проводится так же, как в предыдущей теореме. Здесь для
разнообразия мы определим оператор $L$ не с помощью квадратичной формы,
а как оператор из $\Ho^m(\Om)$ в $H^{-m}(\Om)$. Заметим, что это
определение корректно, т.к. умножение на тест--функцию корректно определено
в соболевских пространствах с негативным индексом гладкости, т.е. в условиях
теоремы на коэффициенты $c_{\al,\be}(x)$ мы имеем $Lu\in H^{-1}(\Om)$ при
$u\in D(\Om)$. Конечно, сужение оператора $L$, действующего в $H^{-1}(\Om)$,
на пространство $\H=L_2(\Om)$ совпадает с оператором, который ставится в
соответствие квадратичной форме $(Lu,u)$ согласно первой теореме о
представлении. Теорема доказана.
\end{proof*nodot}

\subsection{Обобщенная задача Неймана и третья краевая задача}

В предыдущем пункте мы построили оператор $-\De+q(x)$
(и его аналог для сильно эллиптического случая), отвечающий условиям
Дирихле на границе области $\Om$. Мы задаемся вопросом: можно ли построить
оператор $-\De+q$ с сингулярным потенциалом, который отвечает краевым
условиям Неймана или третьей краевой задаче? В этом пункте мы ответим
на этот вопрос.

Рассмотрим сначала одномерный случай: $n=1$, $\Om=(a,b)\subset\R$, который
был изучен в статье Савчука и Шкаликова~\cite{SaSh}. Пусть
$q(x)\in H^{-1}_2(a,b)$, т.е. существует $u(x)\in L_2(a,b)$ такая, что
$u'(x)=q(x)$ в смысле распределений на $(a,b)$. Перепишем выражение
$l(y)=-y''+q(x)y$ в виде
$$
 l(y)=-(y^{[1]})'-u(x)y^{[1]}-u^2(x)y,
$$
где
$$
y^{[1]}=y'-u(x)y
$$
(выражение $y^{[1]}$ называют квазипроизводной).

Рассмотрим линеал $\N$, состоящий из абсолютно непрерывных функций,
для которых $y^{[1]}$ также абсолютно непрерывна. Тогда для $y\in\N$
выражение $l(y)$ корректно определяет функцию из $L_1$, а $y'\in L_2$.
Кроме того, можно интегрировать по частям
\begin{equation}\label{IV:eq:form1}
(l(y),y)=(y^{[1]},y')-(uy,y')-(uy',y)-y^{[1]}\ov{y}|_a^b.
\end{equation}
Если $u(x)$ --- вещественная функция, то выражение~\eqref{IV:eq:form1}
можно переписать в виде
\begin{equation}\label{IV:eq:form2}
(l(y),y)=(y^{[1]},y^{[1]})-(u^2(x)y,y)-y^{[1]}\ov{y}|_a^b.
\end{equation}
Теперь видно, что если положить
\begin{equation}\label{IV:eq:neiman1}
y^{[1]}(a)=y^{[1]}(b)=0
\end{equation}
или
\begin{equation}\label{IV:eq:3prob1}
y^{[1]}(a)=\al y(a),\quad y^{[1]}(b)=\be y(b),
\end{equation}
то выражения~\eqref{IV:eq:form1} или~\eqref{IV:eq:form2}
определят квадратичные формы, которые допускают замыкание в пространстве
$H^1(a,b)$. Утверждение о замыкаемости этих форм следует из элементарных
оценок
\begin{multline*}
 |(u(x)y',y)|+|(u(x)y,y')|+|\al y^2(a)|+|\be y^2(b)|\le\\
\le\ep\|y'\|+(\ep^{-1}\|u\|^2_{L_2}+M)\|y\|^2,
\end{multline*}
где $\ep>0$ можно взять любым, а $M$ зависит только от $\al$ и $\be$.
Из проведеных оценок вытекает: если $y\in\N$, то $y^{[1]}\in L_2(a,b)$,
а также $y'\in L_2(a,b)$. Стало быть, $\N\subset H^1_2(a,b)$. Интересно
отметить следующий факт: если
\begin{equation}\label{IV:eq:N0}
N_0=\{y\in H^1_2(a,b)|l(y)\in L_2(a,b)\},
\end{equation}
то $y^{[1]}$ абсолютно непрерывна, т.е. $\N_0\subset\N$. Действительно,
если выполнено равенство
$$
 -y''+q(x)y=f(x)\in L_2
$$
в смысле распределений, а $y'\in L_2$, то $y^{[1]}\in L_2$, $u^2(x)y\in L_1$;
следовательно, $(y^{[1]})'\in L_1$, что влечет $y^{[1]}\in H_1^1(a,b)$, т.е.
$y^{[1]}$  абсолютно непрерывна. Поэтому для функций из $\N_0$ корректно
определены следы $y^{[1]}(a)$ и $y^{[1]}(b)$. Стало быть, с выражением
$-y''+q(x)y$ мы можем связать оператор $L_U$  с областью
\begin{equation}\label{IV:eq:NU}
\N_U=\{y\in\N| U_j(y)=0,\ j=1,2\},
\end{equation}
где $U_j(y)$ --- линейные формы вида~\eqref{IV:eq:neiman1}
или~\eqref{IV:eq:3prob1} (отметим, что можно рассмотреть и другие
краевые условия, например, квазипериодические: $y(a)=e^{i\th}y(b)$,
$y^{[1]}(a)=e^{i\th}y^{[1]}(b)$, важно лишь, чтобы выражение
$y^{[1]}\ov{y}|_a^b$
можно было выразить через $y(a)$ и $y(b)$).

Из проведенных рассуждений следует, что квадратичная форма $(L_U y,y)$
оператора $L_U$ будет замыкаться в пространстве $H^1_2(a,b)$; более того,
область ее замыкания будет совпадать со всем пространством $H^1_2(a,b)$
(если краевые условия были вида~\eqref{IV:eq:neiman1}
или~\eqref{IV:eq:3prob1}), а линеал $\N_U$ будет плотен в $H^1(a,b)$.

Построенный оператор $L_U$ может быть назван оператором, соответствующим
``обобщенной'' задаче Неймана (в случае условий~\eqref{IV:eq:neiman1})
или третьей краевой задаче (в случае условий~\eqref{IV:eq:3prob1}). Поэтому
для получения
корректного определения оператора $L_U$ остается решить только одну задачу:
{\it доказать плотность линеала $\N_U$ в пространстве $H^1(a,b)$}.
Отметим, что в работе~\cite{SaSh} конструкция операторов $L_U$ осуществляется
по--другому.

Дальнейший наш план состоит в следующем: приведенную конструкцию мы реализуем в
многомерном случае, а для доказательства плотности области определения
мы используем идею аппроксимации операторами с гладкими потенциалами.

Далее будем предполагать, что потенциал $q(x)$ удовлетворяет
условию~\eqref{IV:eq:qin}, т.е. существует вектор--функция
$V(x)=\{V_1(x),\dots,V_n(x)\}$ такая, что
\begin{equation}\label{IV:eq:Vin}
 V_j(x)\in
\left\{
\begin{aligned}
L_2(\Om)&,\quad\text{если $n=1$},\\
L_{2+\ep}(\Om)&,\quad\text{при некотором $\ep>0$, если $n=2$},\\
L_n(\Om)&,\quad\text{если $n\ge 3$},
\end{aligned}
\right.
\end{equation}
причем $q(x)=\div V(x)$. Выражение $\L(f)=(-\De+q(x))f(x)$ представим в виде
\begin{equation}\label{IV:eq:rewrite}
\L(f)=-\div(\nabla f-Vf)-V\cdot(\nabla f-Vf)-V\cdot Vf.
\end{equation}
Очевидно, аналогом квазипроизводной, определенной при $n=1$, является
выражение
$$
 [\nabla]f=\nabla f-Vf.
$$
Ниже для определенности мы рассматриваем случай $n\ge3$ (рассуждения
при $n<3$ аналогичны) и ограничиваемся построением оператора, отвечающего
обобщенной задаче Неймана
\begin{equation}\label{IV:eq:neimanN}
([\nabla]f\cdot \vec n)|_{\pa\Om}=0
\end{equation}
(условие $([\nabla]f\cdot\vec n+\al(x)f)|_{\pa\Om}=0$, где $\al(x)$ --- гладкая
функция на $\pa\Om$, рассматривается точно так же).
Сначала заметим следующее.
\begin{proposition}\label{IV:pr:trace}
Пусть функция $f\in H^1_2(\Om)$ такова, что $\L(f)\in L_2(\Om)$
в смысле распределений на $\Om$. Тогда
$$
 [\nabla]f\in[H^1_{2n/(n+2)}]^n
$$
(здесь степень $n$ означает $n$--ую декартову степень
 соответствующего
пространства). В частности, корректно определен след на границе $\pa\Om$
функции $[\nabla]f$, принадлежащий пространству
$$
 [H^{(n-2)/(2n)}_{2n/(n+2)}(\pa\Om)]^n.
$$
\end{proposition}
\begin{proof}
Если $\L(f)\in L_2(\Om)$, то (как уже было отмечено)
второе и третье слагаемые в~\eqref{IV:eq:rewrite} принадлежат пространству
$L_{2n/(n+2)}(\Om)\supset L_2(\Om)$. Но тогда слагаемое
$\div[\nabla]f\in L_{2n/(n+2)}(\Om)$, и, следовательно,
$[\nabla]f\in [H^1_{2n/(n+2)}(\Om)]^n$. Утверждение
о следе этой функции теперь
вытекает из известной теоремы о следах (см., например,~\cite{Tr}).
Предложение доказано.
\end{proof}

Теперь определим оператор $L_N$ (здесь индекс $N$ означает, что мы
рассматриваем краевые условия Неймана) следующим образом:
\begin{equation}\label{IV:eq:defLN}
\left\{
\begin{aligned}
L_Nf&=\L(f),
\\
\D(L_N)&=\{f\in H^1(\Om), \L(f)\in L_2, ([\nabla]f\cdot\vec n)|_{\pa\Om}=0\}.
\end{aligned}
\right.
\end{equation}

Из предложения~\ref{IV:pr:trace} следует, что это определение корректно.
В силу формулы Грина получаем, что квадратичная форма оператора $L_N$ равна
$$
 (L_Nf,f)=(\nabla f,\nabla f)-(\ov{U}f,\nabla f)-(\nabla f, Uf).
$$
Проводя элементарные оценки (такие же, как при $n=1$), получим, что эта форма
замыкаема, причем область определения замыкания есть подпространство в
$H^1(\Om)$. Из классической теории известно, что если потенциал $q(x)$ ---
гладкая функция, то область определения замыкания формы совпадает с
$H^1(\Om)$.

Обозначим теперь $\H^1=H^1(\Om)$, а $\H^{-1}$ --- дуальное к $H^1(\Om)$
по отношению к скалярному произведению в $L_2(\Om)$. При $n=1$ дуальное
пространство $\H^{-1}$ описать легко (оно на $2$ размерности шире,
нежели $H^{-1}(\Om)$); при $n\ge2$ описание $\H^{-1}$ нам неизвестно.
Пусть $V(x)\in (L_n(\Om))^n$, т.е. выполнено условие~\eqref{IV:eq:Vin},
и пусть $V_k(x)$ --- последовательность
функций из $(E(\Om))^n$, таких, что $V_k(x)$ сходятся к $V(x)$ в
норме $(L_n(\Om))^n$. Положим $q_k=\div V_k$, $q=\div V$, причем
в последнем случае $\div$ понимается в смысле распределений из $D'(\Om)$.
Тогда $q_k(x)$ --- последовательность функций из $E(\Om)$, сходящихся
к $q(x)$ в норме $H^{-1}_n$. Обозначим через $L_k$ операторы $-\De+q_k$,
порождаемые в пространстве $L_2(\Om)$ краевыми условиями
$$
 (\nabla f-V_k f,\vec n)|_{\pa\Om}=0.
$$
Квадратичные формы $(L_k u,u)=(\nabla u,\nabla u)-(\nabla u,V_k u)
-(\ov{V_k} u,\nabla u)$, определенные на $\D(L_k)$, секториальны,
причем области определения их замыканий $l_k[u]$ совпадают с $\H^1$.
Введем также секториальную форму
$$
  l[u]=(\nabla u,\nabla u)-(\nabla u,Vu)-(\ov Vu,\nabla u),
$$
определенную на $\H^1$. Поскольку $V\in L_n(\Om)$ и $V_k\to V$ в норме
этого пространства, то имеет место сходимость форм $l_k[u]\to l[u]$ при всех
$u\in H_2^1(\Om)$.

Пусть $\hat L$ --- $m$--секториальный оператор в $\H$, ассоциированный
с формой $l[u]$ согласно первой теореме о представлении.
Тогда $\hat L+\rho$ ограниченно обратим (при достаточно больших $\rho$),
и гомеоморфно отображает $\D(\hat L)$ на $L_2$. В дальнейшем, не
ограничивая общности, полагаем $\hat L$, а также все операторы $L_k$,
равномерно ограниченно обратимыми.

Пусть $T_0=-\De+1$ --- самосопряженный равномерно положительный оператор в
$L_2$, порожденный краевым условием $(\nabla f,\vec n)|_{\pa\Om}=0$.
Тогда $T_0$ гомеоморфно переводит $\H^1$ в $\H^{-1}$, а $T_0^{1/2}$
также гомеоморфно отображает $\H^1$ в $\H^0$ и $\H^0$ в $\H^{-1}$.
При этом оператор $\hat L$ допускает представление
$$
\hat L=T_0^{1/2} (A+iB) T_0^{1/2},
$$
где $A$ --- положительно определенный ограниченный, а $B$ --- самосопряженный
ограниченный оператор в $\H^0$. Теперь, по аналогии с
теоремой~3 статьи~\cite{NeSh}, получаем, что операторы $L_k$ сходятся
к $\hat L$ в смысле равномерной резольвентной сходимости.

Докажем, что оператор $\hat L$ и оператор $L_N$, определенный
в~\eqref{IV:eq:defLN}, совпадают.

Зафиксируем $g\in L_2(\Om)$ и положим $f=\hat L^{-1}g$.
По квадратичной форме оператора $\hat L$ однозначно восстанавливается
билинейная форма на $\H^1\times\H^1$, поэтому
для произвольной функции $\phi\in D(\Om)$ выполнено $(\hat Lf,\phi)=
(\nabla f,\nabla \phi)-(\nabla f, V\phi)-(\ov V f,\nabla\phi)$.
С другой стороны, $\L(f)$ для $f\in\H^1$ определено в смысле
распределений на $\Om$, и для $f\in\H^1$ и $\phi\in D(\Om)$
по формуле Грина также имеем
$$
(L(f),\phi)=(\nabla f,\nabla \phi)-(\nabla f, V\phi)-(\ov V f,\nabla\phi)=
(g,\phi).
$$
Поэтому $L(f)\in L_2(\Om)$. Из предложения~\ref{IV:pr:trace} следует, что
обобщенный градиент $[\nabla]f$ принадлежит пространству
$H^1_{2n/(n+2)}$ и имеет след на границе.

Пусть теперь $f_k=L_k^{-1}g$. Из резольвентной сходимости следует, что
$f_k\to f$ в норме $\H^1$. Докажем, что $\nabla f_k-V_k f_k\to [\nabla] f$
в норме $H^1_{2n/(n+2)}$. Тогда из того, что
$(\nabla f_k-V_k f_k)\cdot \vec n=0$,
будет следовать, что $(\nabla f-V f)\cdot \vec n=0$, т.е. $f\in\D(L_N)$.

Поскольку $L_kf_k=g$ и $L_Nf=g$, то
\begin{multline*}
0=L_N f-L_k f_k=-\div(\nabla f-Vf-(\nabla f_k-V_k f_k))-
V(\nabla f-Vf-\\
-(\nabla f_k-V_k f_k))-(V-V_k)(\nabla f_k-V_k f_k)-
(V^2-V_k^2)f.
\end{multline*}
Непосредственно проверяется, что в последнем выражении все слагаемые,
кроме первого, сходятся в норме пространства $L_{2n/(n+2)}$. Тогда
и первое слагаемое сходится в этой же норме, а поэтому
выражение $(\nabla f-Vf-(\nabla f_k-V_k f_k))$, находящееся под $\div$,
сходится в норме $H^1_{2n/(n+2)}$.

Итак, $f\in\D(L_N)$. Если $g$ пробегает все пространство
$L_2(\Om)$, то $f$ пробегает все $\D(\hat L)$. Поэтому
$\D(\hat L)\subset \D(L_N)$. Из секториальности оператора $L_N$ и
$m$--секториальности $\hat L$ теперь следует,
что $L=\hat L$.

Таким образом, мы доказали следующий результат:
\begin{theorem}\label{IV:th:mainNeim}
Пусть функция $q(x)$ такова, что выполнено
условие~\eqref{IV:eq:qin}.
Тогда оператор $L_N$ корректно определен формулой~\eqref{IV:eq:defLN}, причем
его область определения $\D(L_N)$ плотна в $H^1_2(\Om)$. В случае
вещественности функции $q(x)$ оператор $L_N$ является самосопряженным
и полуограниченным, причем при достаточно больших $\rho$ область
$\D((L_N+\rho)^{1/2})$ совпадает c $H^1_2(\Om)$.
\end{theorem}

\subsection{Асимптотика собственных значений операторов $L_D$ и $L_N$}

Пусть $0\le p\le 1$, а $T$ --- равномерно положительный самосопряженный
оператор в $\H$. Оператор $Q$ называют $p$--подчиненным оператору $T$,
если $\D(Q)\subset\D(T^p)$, причем
$$
 |(Qx,x)|\le C(T^p x,x)\quad\forall x\in\D(T^p),
$$
где постоянная $C$ не зависит от $x$.

Для $p$--подчиненных операторов известны теоремы о сохранении главных
членов асимптотики собственных значений и о свойствах базисности
собственных функций. Теория подчиненных операторов развивалась
А.~С.~Маркусом, В.~И.~Мацаевым, В.~Э.~Кацнельсоном, М.~С.~Аграновичем,
А.~А.~Шкаликовым и другими.

Положим $T_0=-\De$ в пространстве $\H=L_2(\Om)$ ($\Om$ --- ограниченная
гладкая область в $\R^n$) с областью $\D(T_0)=H^2(\Om)\cap\Ho^1(\Om)$.
С оператором $T_0$ связываем шкалу гильбертовых пространств $\H^\th$,
$\th\in\R$, причем $\H^1=\Ho^1(\Om)=\D(T_0^{1/2})$, $\H^{-1}=H^{-1}(\Om)$.

\begin{lemm}
Пусть $q(x)$ есть мультипликатор из пространства $\H^\th$ в пространство
$\H^{-\th}$
при некотором $0\le\th<1$ (будем писать $q(x)\in M[\th,\Om]$).
Тогда оператор ``умножения'' на функцию $q(x)$ является $\th$--подчиненным
оператору $T_0$.
\end{lemm}
\begin{proof}
Сформулированное утверждение является лишь перефразировкой понятий.
Как и ранее (но для случая пространства $\R^n$), мы считаем, что
$q(x)\in M[\th,\Om]$, если
$$
  |(qf(x),g(x))|\le C\|f\|_\th\|g\|_\th\quad f,g\in D(\Om).
$$
Полагая в этом неравенстве $f=g$ и принимая во внимание, что
$\|f\|_\th=\|T^{\th/2}f\|$ (согласно определению шкалы пространств),
получаем
$$
  |(q f(x), f(x))|\le C\|T^{\th/2} f\|^2=C(T^\th f,f),\quad f\in D(\Om).
$$
Заметим, что $D(\Om)$ плотно в $\H^1$, и, тем более, в $\H^\th$ при $\th<1$.
Поэтому последнее неравенство распространяется на $f\in\H^\th$, тем
более, на $f\in \D(T^\th)$. Лемма доказана.
\end{proof}

\begin{theorem}\label{IV:th:Ndistr}
Пусть $q\in M[\th,\Om]$ при некотором $\th<1$. Пусть $N(\la,L_D)$ ---
функция распределения собственных значений определенного в п.1 оператора
$L_D=-\De+q(x)$. Тогда при $r\to\infty$ справедлива асимптотика
\begin{equation}\label{IV:eq:Ndistr}
N(r,L_D)=(2\pi)^{-n}\om r^{n/2}+O(r^{n/2+\th-1})+O(r^{(n-1)/2}),
\end{equation}
где $\om$ --- объем области $\Om\subset\R^n$.
В частности, при $n\ge2$ эта формула справедлива, если
$q(x)\in H^{-\th}_p(\Om)$ при $p\ge n/\th$.

В случае $n=1$ формула~\eqref{IV:eq:Ndistr} справедлива при
\begin{equation}\label{IV:eq:qcond}
q(x)\in\left\{
\begin{array}{ll}
H^{-\th}_2(\Om),&\quad\text{если $1/2<\th<1$},\\
H^{-\th}_p(\Om),\ p\ge n/\th,&\quad\text{если $0<\th<1/2$},\\
H^{-1/2}_p(\Om),\ p>2,&\quad\text{если $\th=1/2$}.
\end{array}
\right.
\end{equation}
Наконец, если функция $q(x)$ такова, что
\begin{equation}\label{IV:eq:qcond1}
q(x)\in\left\{
\begin{array}{ll}
H^{-1}_2(\Om)&,\quad\text{если $n=1$},\\
H^{-1}_p(\Om),\ p>2,&\quad\text{если $n=2$},\\
H^{-1}_n(\Om)&,\quad\text{если $n\ge 3$},
\end{array}
\right.
\end{equation}
то
\begin{equation}\label{IV:eq:Ndistrweak}
N(r,L_D)=(2\pi)^{-n}\om r^{n/2}(1+o(1)),\quad r\to\infty.
\end{equation}
\end{theorem}
\begin{proof}
Функция распределения собственных значений оператора $T_0=-\De$,
порожденного краевым условием Дирихле, имеет асимптотику
(см. например,~\cite{Vas})
\begin{equation}\label{IV:eq:Ndistrmodel}
N(r,T_0)=(2\pi)^{-n}\om r^{n/2}+O(r^{(n-1)/2})),\quad r\to\infty.
\end{equation}
Более того, известно~\cite[IV]{Vas}, что при некоторых дополнительных
условиях на область $\Om$ остаток $O(r^{(n-1)/2})$ может быть уточнен:
он равен $h r^{(n-1)/2}(1+o(1))$, где $h=\const$.

Из леммы 5 и теоремы Маркуса--Мацаева~\cite{MM1} (см. также~\cite{Ag1})
получаем, что для функции $N(r,L_D)$ возмущенного оператора $L_D=T_0+q$
справедлива оценка
$$
 |N(r,T_0)-N(r,L_D)|\le C(N(r+cr^\th,T_0)-N(r-cr^\th,T_0)),
$$
где $C$ и $c$ --- некоторые константы. Учитывая~\eqref{IV:eq:Ndistrmodel},
получаем, что правая часть последнего неравенства допускает оценку
$\le C(r^{n/2+\th-1})+Cr^{(n-1)/2}$. Это влечет оценку~\eqref{IV:eq:Ndistr}.

Второе утверждение теоремы (о достаточном условии принадлежности
$q(x)\in M[\th,\Om]$) вытекает из результатов работы~\cite{BSh}
(как уже говорилось в п.1, результаты, полученные при
$\Om=\R^n$, переносятся на случай ограниченной гладкой области $\Om$,
если мы имеем дело с задачей Дирихле).

Наконец, если выполняется условие~\eqref{IV:eq:qcond1} (это отвечает
случаю $\th=1$), то оператор $L_D=-\De+q(x)$ является относительно
компактным возмущением $T_0$, а тогда формула~\eqref{IV:eq:Ndistrweak}
следует из теоремы М.~В.~Келдыша (см.~\cite[гл.1]{Ma}).
\end{proof}

\begin{corollary}
Если условие теоремы~\ref{IV:th:Ndistr} выполнено при
$\th\le 1/2$, то
\begin{equation}\label{IV:eq:Ndistrclassic}
N(r,L_D)=(2\pi)^{-1}\om r^{n/2}+O(r^{(n-1)/2}),
\end{equation}
т.е. остаток остается таким же, как для оператора Лапласа с условием
Дирихле.
\end{corollary}
\begin{proof*nodot} вытекает из~\eqref{IV:eq:Ndistr} и неравенства
$n/2+\th-1\le (n-1)/2$, если $\th\ge1/2$.
\end{proof*nodot}

Теперь займемся вопросом о распределении собственных значений
построенного в п.2 оператора $L_N$, отвечающего обобщенной
задаче Неймана.

Обозначим через $T_N$ оператор Лапласа с краевым условием Неймана,
и через $L_N$ --- его возмущение потенциалом--распределением $q(x)$.
Имеем
$$
  (L_N u,u)=\|\nabla u\|^2-(\nabla u,V u)-(\ov V u,\nabla u),\quad
u\in\D(L_N).
$$
Если потенциал $q$ таков, что выполнены условия
теоремы~\ref{IV:th:mainNeim}, то
$$
 |(\nabla u, Vu)|\le \ep\|\nabla u\|^2+C\|u\|^2,
$$
а потому $L_N$ есть компактное возмущение оператора $T_N$. Поскольку
главный член асимптотики собственных значений и оценка остатка
для оператора $T_N$ такие же, как для оператора $T_0$, то
в условиях теоремы~\ref{IV:th:mainNeim} с помощью теоремы Келдыша
получаем
\begin{equation}\label{IV:eq:NdistrweakL}
N(r,L_N)=(2\pi)^{-1}\om r^{n/2}(1+o(1)),\quad r\to\infty.
\end{equation}

Для получения аналога утверждения теоремы~\ref{IV:th:Ndistr} нужно найти
условия на функцию $V(x)$, при некоторых справедлива оценка
\begin{equation}\label{IV:eq:Vestim}
|(\nabla u, V u)|\le C\|u\|^2_\th,\quad\th<1,
\end{equation}
где $C$ зависит только от $V$ и $\th$, а $\|\cdot\|_\th$ --- норма
в $H_2^\th(\Om)$. Очевидно, именно такая оценка гарантирует
$\th$--подчиненность оператора $L_N$ оператору $T_N$ в смысле
квадратичных форм.

Конечно, несложно получить грубые достаточные условия на функцию $V(x)$,
при которых выполняется оценка~\eqref{IV:eq:Vestim}. Но наша гипотеза
состоит в том, что эта оценка выполняется в условиях
теоремы~\ref{IV:th:mainNeim}. Здесь мы приведем доказательство этой
гипотезы только при $\th<1/2$.

Действительно, согласно теореме Гривара (см.~\cite{Grivar},
и~\cite[гл.5]{Tr}) имеем, что при $\th<1/2$ совпадают интерполяционные
пространства
$$
 [\Ho_2^1(\Om),L_2(\Om)]_\th=[H^1_2(\Om),L_2(\Om)]_\th=:H_2^\th(\Om).
$$
Поскольку пространство основных функций $D(\Om)$ плотно в
$\Ho_2^1(\Om)$, то $D(\Om)$ также плотно в $H^\th_2(\Om)$ при $\th<1/2$.
Следовательно, в этом случае достаточно доказать
оценку~\eqref{IV:eq:Vestim} для $u\in D(\Om)$. Но для таких функций
$u$ имеем $(\nabla u, V u)=(u,q u)$, где $q=\div V$. Так как
нормы $\|T_0^{\th/2}u\|$ и $\|u\|_\th$ при $\th<1/2$ эквивалентны,
то оценка~\eqref{IV:eq:Vestim} вытекает из леммы 5.
Тем самым мы доказали следующий результат:
\begin{theorem}\label{IV:th:NdistNeim}
Пусть функция $q(x)$ подчинена условию~\eqref{IV:eq:qcond1}. Тогда для
функции распределения собственных значений $N(r,L_N)$ выполнено
соотношение~\eqref{IV:eq:NdistrweakL}. Если же $q(x)\in H_p^{-\th}(\Om)$
при $p\ge n/\th$, $\th<1/2$, то справедлива
формула~\eqref{IV:eq:Ndistrclassic}.
\end{theorem}

\subsection{Свойства собственных функций операторов $L_D$ и $L_N$}
Если функция $q(x)$ вещественна, то построенные операторы $L_D$ и $L_N$
являются самосопряженными, а потому их собственные функции, согласно теореме
Гильберта--Шмидта, образует ортогональный базис пространства $L_2(\Om)$.
Конечно, это не так, если $q(x)$ --- комплексный потенциал. Однако
при выполнении условия~\eqref{IV:eq:qcond1} операторы $L_D$
и $L_N$ являются относительно компактными возмущениями
самосопряженных операторов $T_0$ и $T_N$, а потому система
их собственных и присоединенных функций, согласно теореме Келдыша
(см.~\cite{Ke},\cite[гл.III]{GK}), будет образовывать полную систему
(здесь мы учитываем, что операторы $T_0^{-1}$ и $T_N^{-1}$
являются операторами конечного порядка, т.е. операторы $T_0^{-k}$
и $T_N^{-k}$ являются ядерными при достаточно большом $k$). Для
$\th$--подчиненных операторов при $\th<1$ можно утверждать большее.

Далее для простоты изложения будем считать, что присоединенных функций у
операторов $L_D$ и $L_N$ нет. Напомним определение базиса для
суммирования методом Абеля, принадлежащий В.~Б.~Лидскому~\cite{Lid}.
Пусть $\{y_k\}_1^\infty$ --- система собственных функций оператора $A$,
действующего в пространстве $\H$, а $\{z_k\}_1^\infty$ ---
биортогональная система (согласно теореме Келдыша~\cite{Ke}, такая система
существует для любого оператора $A$ с дискретным спектром). Система
$\{y_k\}$ образует базис для метода суммирования Абеля порядка $\al>0$,
если ряд
$$
 f(t)=\suml_{k=1}^\infty(f,Z_k)e^{-\la_k^\al t}y_k
$$
сходится при всех $t>0$ и $f\in\H$ (после возможной расстановки скобок,
не зависящей от $t$ и $f$), причем $\|f(t)-f\|\to0$ при $t\to+0$. Здесь
выбирается главная ветвь функции $\la^\al$ (положительная при $\la>0$)
с разрезом по отрицательной оси.

Сформулируем теперь результат о полноте и базисности собственных функций.

\begin{theorem}\label{IV:th:Evec}
Пусть $q(x)$ такова, что выполнено условие~\eqref{IV:eq:qcond1}. Тогда
системы собственных и присоединенных функций операторов $L_D$ и $L_N$
образуют полную систему в $L_2(\Om)$. Более того, эти системы образуют
базис для метода суммирования Абеля порядка $\al\ge n/2$. При дополнительных
предположениях порядок суммирования можно понизить. А именно, если
$n\ge 2$, $\th<1$, выполнено условие $q(x)\in H^{-\th}_p(\Om)$, $p\ge n/\th$
или условие~\eqref{IV:eq:qcond} при $n=1$, то система собственных и
присоединенных функций оператора $L_D$ образует базис для метода
суммирования Абеля любого порядка $\al>\dfrac{n}2-(1-\th)$. В случае
$\th<1/2$ такое же утверждение справедливо для системы собственных и
присоединенных функций оператора $L_N$.
\end{theorem}

\begin{proof*nodot} первого утверждения о полноте уже было проведено.
Доказательство второго утверждения вытекает из теоремы Лидского~\cite{Lid}
с учетом дополнений к этой теореме, сделанных Мацаевым (см.~\cite{Ag2}),
и Шкаликовым~\cite{Shkal}. Наконец, последнее утверждение следует из
теорем Маркуса--Мацаева--Аграновича (см.~\cite{Ag1}).
\end{proof*nodot}

\end{document}